\documentclass [12pt,a4paper] {amsart}

\usepackage[utf8]{inputenc}
\usepackage[T1]{fontenc}
\usepackage{lmodern}
\usepackage{microtype}

\usepackage {amssymb}
\usepackage {amsmath}
\usepackage {amsthm}
\usepackage{mathtools}
\usepackage{xfrac}

\usepackage{url}
 
\usepackage[titletoc,toc,title]{appendix}

\usepackage{enumitem}
\usepackage{hyperref}
\hypersetup{
    colorlinks=true,
    linkcolor=blue,
    filecolor=magenta,      
    urlcolor=cyan,
    citecolor=magenta
}
\usepackage{cleveref}

\usepackage{tikz}
\usetikzlibrary{arrows}

\usepackage[OT6,T1]{fontenc}
\usepackage[utf8]{inputenc}

\newcommand{\armenian}{\fontencoding{OT6}\fontfamily{cmr}\selectfont}
\DeclareTextFontCommand{\textarmenian}{\armenian}

\theoremstyle {plain}

\theoremstyle {definition}

\theoremstyle {remark}

\usepackage{geometry}

\calclayout

\title{The mathematical work of Anania Shirakatsi}
\author{Vahagn Aslanyan}

\email{Vahagn.Aslanyan@manchester.ac.uk}
\address{Department of Mathematics, University of Manchester, Manchester, UK}

\begin{document}

\vspace*{-1cm}

\thanks{This work was supported by EPSRC Open Fellowship EP/X009823/1 and Dame Kathleen Ollerenshaw Fellowship.  For the purpose of open access, the author has applied a Creative Commons Attribution (CC BY) licence to any Author Accepted Manuscript version arising from this submission.}


\subjclass[2020]{}

\maketitle

\begin{abstract}
This paper aims to give a brief account of the mathematical work of the 7th-century Armenian polymath and natural philosopher Anania Shirakatsi. The three sections of Anania's ``Book of Arithmetic'' -- tables of arithmetic operations, a list of problems and their answers, and a collection of entertaining puzzles -- are presented and discussed, the focus being on the problems and solutions. A close examination of the structure of these problems reveals that Anania's proficiency in arithmetic was considerably more sophisticated than their mathematical contents might suggest. The geography of Anania's problems is illustrated through two maps highlighting the locations referenced within these problems.
\end{abstract}

\section{Introduction}

In his book \textit{History of the Armenians} Movses Khorenatsi, a 5th-century Armenian historian (considered to be the father of Armenian history), writes, ``For even though we are small and very limited in numbers and have been conquered many times by foreign kingdoms, yet too, many acts of bravery have been performed in our land, worthy of being written and remembered.''\footnote{The English version of this quote is taken from the Wikipedia page \href{https://en.wikipedia.org/wiki/Movses_Khorenatsi}{Movses Khorenatsi}. The original version is as follows: ``\textarmenian{Զի թէպէտ եւ եմք ածու փոքր եւ թուով յոյժ ընդ փոքու սահմանեալ, եւ զօրութեամբ տկար, եւ ընդ այլով յոլով անգամ նուաճեալ թագաւորութեամբ՝ սակայն բազում գործք արութեան գտանին գործեալ եւ ի մերում աշխարհիս, եւ արժանի գրոյ յիշատակի}:''} Similarly, I think that even though mathematics in Armenia in ancient and medieval times was not as advanced as, say, in Babylon, Greece, Egypt, India, Persia, and China, there have been important mathematical developments in Armenia too that are worth studying and writing about. This paper aims to contribute to this goal by presenting and analysing the mathematical work of one of the earliest and most prominent Armenian mathematicians, Anania Shirakatsi.

Anania Shirakatsi (c. 610 -- c. 685), also known as Ananias of Shirak in the English-speaking world, was a 7th-century Armenian polymath who is considered the father of science in Armenia and the first Armenian mathematician. His work covered many areas including mathematics, geography, astronomy, and chronology. His most famous mathematical work is the \textit{Book of Arithmetic} which consists of three sections and which is the subject of study in this paper. The first section contains tables of addition, multiplication, subtraction, and division, and some short explanations of these arithmetic operations. The second section consists of 24 problems along with their answers. The third section contains entertaining mathematical puzzles intended to be used during social gatherings to amuse and challenge one another.

Anania's mathematical work has been studied by several Armenian scholars, including the historian Abrahamyan, historian of mathematics Petrosyan, and Armenologist Matevosyan \cite{abrahamyan-shirakaci-matenagrutyun,abrahamyan-petrosyan-shirakaci-matenagrutyun,matevosyan-bazhanum,matevosyan-norahayt-patarikner}. Several decades before them the orientalist/historian Orbeli translated Anania's problems and solutions from the Book of Arithmetic into Russian and published both the Armenian and Russian versions in one book \cite{orbeli-shirakatsi}. There are several papers and book chapters in English about Anania's mathematical work and its historical context. For instance, Shaw \cite{shaw-overlooked-numeral-system} and Chrisomalis \cite{Chrisomalis-numerical-notation} wrote about the Armenian numeral systems including the one developed and used by Anania, Matthews \cite{matthews-anania-of-shirak} wrote a brief review of Anania's work, Greenwood \cite{greenwood-reassessment} provided a thorough overview of Anania's life, as well as his problems and solutions and their historical context, and Hewsen \cite{hewsen-science-7th-cent-arm} wrote on Anania's scientific texts briefly touching upon his mathematical work.

In spite of that, there seems to be no study in English of all three parts of the Book of Arithmetic. This paper aims to bridge this gap by providing a detailed overview of the Book of Arithmetic in its entirety. However, I mainly focus on Anania's problems and solutions, and offer an analysis of their mathematical, historical, and geographical contents. The problems and solutions were translated into English by Greenwood \cite{greenwood-reassessment}, and I am not aware of any other English text containing all of the problems. Greenwood used Abrahamyan's text \cite{abrahamyan-shirakaci-matenagrutyun} as a source, while my translations are based on the books of Orbeli \cite{orbeli-shirakatsi} and Abrahamyan-Petrosyan \cite{abrahamyan-petrosyan-shirakaci-matenagrutyun}. 

Anania's work has survived and reached us through various manuscripts copied from the original or from earlier copies which are kept in Matenadaran (the Mesrop Mashtots Institute of Ancient Manuscripts)\footnote{It is possible to take a \href{https://matenadaran3d.am/en}{virtual tour} of the Matenadaran museum where (copies of) Anania's works can be seen among numerous other invaluable manuscripts.}; see Figure \ref{fig:Arithmetic} for a page of a 13th-century manuscript containing tables of arithmetic operations from Anania's book. Some of the above-mentioned books and papers used these manuscript as their sources. While different manuscripts disagree at places, they are mostly identical, and contain several minor inaccuracies in the statements of Anania's problems. Greenwood gave literal translations and rightly noted the unsolvability of some of the problems (7 in number). However, in most cases it is clear how the statements of problems are supposed to be interpreted from a mathematical point of view, hence to a mathematician these problems make perfect sense and are solvable. For other problems, since Anania  supplied solutions, or rather answers, to these problems, it allows one to fix the inaccuracies in the statements and again get complete and meaningful problems. This was mostly done in \cite{abrahamyan-petrosyan-shirakaci-matenagrutyun}, which is one of the sources I used for translating the problems. My translations are not literal, but they are mathematically coherent and still reflect the historical and geographical contents of the problems. This is indeed important; as Greenwood \cite{greenwood-reassessment} points out, Anania's problems are ``a rich source for seventh-century history whose value has not been sufficiently recognised''. With some of the problems I provide brief comments on the historical events mentioned in them, as well as point out the changes made to make the problems mathematically complete. After presenting the 24 problems that are commonly considered to be part of Anania's Book of Arithmetic, I add three further problems translated from \cite{matevosyan-norahayt-patarikner} which are arguably due to Anania too.

Further, I provide two maps highlighting the geographic locations mentioned in Anania's problems; see Figures \ref{fig:Sasanian Armenia} and \ref{fig:Ayrarat map}.  The first is a map of Sasanian (Persian) Armenia in the 5th-6th century, when the historical events of Anania's problems took place. The second map shows the historical Ayrarat province of Greater Armenia, which later became part of Sasanian Armenia. Most of the locations in Anania's problems are indeed in Ayrarat and are all marked on this map.

I also analyse the mathematical contents of the problems and argue that Anania's arithmetical skills and knowledge were more advanced than what can directly be seen in the Book of Arithmetic.

The paper is organised as follows. In Section \ref{sec: tables} I give a brief overview of the tables of arithmetic operations found in the Book of Arithmetic. Section \ref{sec: prob and sol} contains Anania's problems and solutions, their mathematical analysis, and the two above-mentioned maps. In Section \ref{sec: entertaining puzzles} I present translations of Anania's entertaining puzzles.

\section{Tables of arithmetic operations}\label{sec: tables}

\begin{figure}
    \centering
    \includegraphics[width=1\linewidth]{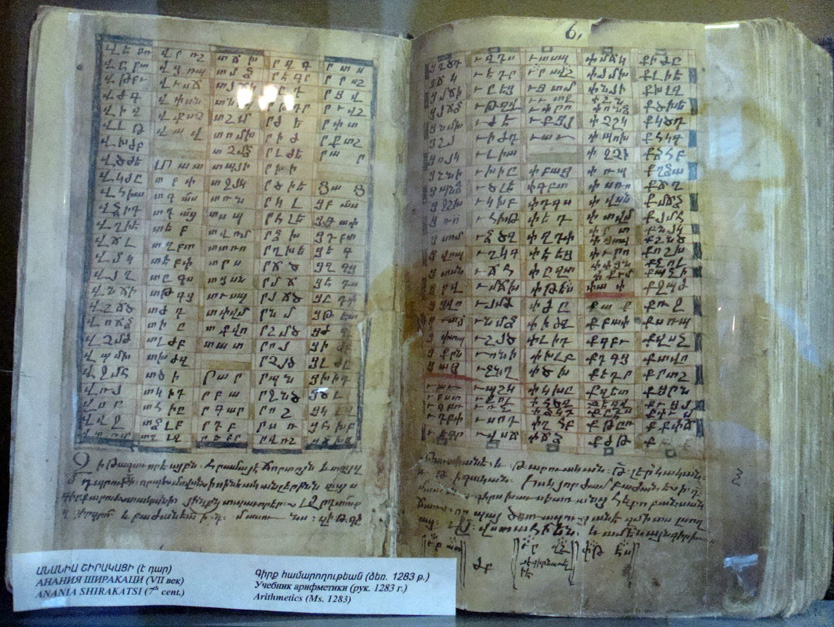}
    \caption{Anania's \textit{Book of Arithmetic}. Manuscript dated to 1283, currently kept in Matenadaran in Yerevan.\\ \scriptsize{The image was downloaded from Wikipedia and is \href{https://commons.wikimedia.org/wiki/File:Shirakatsi_manuscript.JPG}{copyright free}. The original uploader was Taron Saharyan at Russian Wikipedia, Public domain, via Wikimedia Commons.}}
    \label{fig:Arithmetic}
\end{figure}

The \textit{Book of Arithmetic} contains a number of tables of the four basic arithmetic operations -- addition, multiplication, subtraction, and division. An example of such a table can be seen in Figure~\ref{fig:Arithmetic}. It is worth pointing out that, as can be seen in the picture, Anania used the letters of the Armenian alphabet to denote numbers. However, this was different from the traditional representation of numbers in Armenian texts which was analogous to the Greek numeral system, i.e. the first nine letters of the alphabet correspond to the numbers 1--9, the next nine letters correspond to 10--90, and so on (e.g. \textarmenian{իգ} is 23).\footnote{The Armenian system enjoyed the luxury of having 36 letters, as opposed to Greek 24, which allowed one to write numbers up to 9999 without introducing any additional symbols.} The unique numeral system that Anania developed was based on 12 letters. It used the first nine letters of the alphabet (\textarmenian{ա,բ,գ,դ,ե,զ,է,ը,թ}) for the units and the letters \textarmenian{ժ, ճ, ռ} to denote 10, 100, 1000 respectively (as in the standard system). Numbers up to 9,999 were written using these 12 letters and a multiplicative-additive notation, e.g. 50 was written as \textarmenian{եժ}, and 216 was written as \textarmenian{բճժզ}. Using the other letters of the alphabet, Anania was able to represent all numbers up to 999,999. Larger numbers were written using special symbols such as $\hat{}$, for instance, $\hat{\textarmenian{Ք}}=9\cdot 10^7$. Often (e.g. in the problems and answers, see \ref{sec: prob and sol}) numbers were written in words.

The tables of addition list all sums of the form $a\cdot 10^n + b\cdot 10^n$ where $1\leq a \leq b \leq 9,~ 0 \leq n \leq 3$. The tables of multiplication list all products of the form $(a\cdot 10^m) \cdot (b\cdot 10^n)$ where $1\leq a, b \leq 9,~ 0 \leq m, n \leq 3$, and in some tables $b$ also takes the value $10$. The tables of subtraction list all differences which can be obtained from the tables of addition. The tables of division list the ratios of the form $\frac{6000}{n}$ where $1\leq n \leq 100$ or $n=a\cdot 100$ or $n= a\cdot 1000$ with $1\leq a \leq 9$. When the ratio is not an integer, the closest integral value is recorded, e.g. $353$ for $\frac{6000}{17}$. 

The largest number occurring in all these tables is $9\cdot 10^{7}$.

There are also tables containing even and odd numbers. The numbers $a\cdot 10^n$ with $a=2,4,6,8$ and $0\leq n \leq 10$ are listed as even, and the numbers $a\cdot 10^n$ with $a=1,3,5,7,9$ and $0\leq n \leq 10$ are listed as odd\footnote{Clearly, Anania had a slightly odd notion of an odd number.}.

\section{Problems and solutions}\label{sec: prob and sol}

In this section I present the problems and solutions of Anania Shirakatsi together with some commentary. Some of the problems (e.g. the 23rd one) are incomplete and cannot be solved, but the solutions can be used to recover the missing information in them.

\subsection{Problems}

\begin{enumerate}[itemsep=1ex]
    \item My father told me that during the war against the Persians,  Zorak Kamsarakan was distinguished for his bravery. Within a month he attacked the Persian troops three times. During the first attack he destroyed half of the enemy troops, the second time he killed a quarter, and the third time -- an eleventh. The rest, two hundred eighty in number, fled to Nakhchavan. Now, we must work out the number of Persian troops before the attacks. \label{prob:1}

    \noindent \textbf{Comments.} \textit{Kamsarakans were an Armenian noble family that ruled over two cantons of the historical Ayrarat province of Armenia --  Shirak and Arsharunik. Their capital was the city of Yervandashat. Zorak Kamsarakan was a member of this family who lived in the 6th century when Armenia was divided between the Byzantine and Sasanian Empires, with Ayrarat  under the Sasanian rule (until 591). He took part in the rebellion of 571-572 against the Persians. He is known to have attacked the Persian troops three times causing significant losses. In particular, during the battle of Dvin on 2nd February 572, Zorak destroyed the forces (15,000 strong) of Suren, the Marzpan (margrave/governor) of Sasanian Armenia, then killed Suren, and liberated Dvin. Nakhchavan was one of the cantons of the province of Vaspurakan (there was also a city of the same name). See Figures \ref{fig:Sasanian Armenia} and \ref{fig:Ayrarat map}.}

    \item One of my friends travelled to Bahl and bought pearls at a profitable price. When he returned home in Ganjak, he sold half of the pearls for 50 drams each. Then he went to Nakhchavan and sold a quarter of the pearls for 70 drams each. After that he travelled to Dvin and sold a twelfth of the pearls for 50 drams each. When he visited us in Shirak, he had 24 pearls left. Now, using this information work out the total number of pearls and how much money my friend made selling them.

     \noindent \textbf{Comments.} \textit{The dram was a unit of mass and currency. The current currency of Armenia is the Armenian dram.}

  \textit{Bahl refers to Balkh, a city in Sasanian Empire, in present-day Afghanistan. In ancient times it was known as Bactra and was the capital of Bactria (Greco-Bactrian kingdom). Ganjak (or Ganzak) refers to Ganjak Sahastan, a historical city in Sasanian Empire that used to be the capital of Media Atropatene. It was situated south-east of lake Urmia, in present-day Iran. Dvin was a historical capital of Armenia (336-428 AD) in the Ayrarat province. It is in present-day Armenia. See Figures \ref{fig:Sasanian Armenia} and \ref{fig:Ayrarat map}.}

    \item My teacher told me that thieves robbed the treasury of Marcian and stole a half and a quarter of the treasure. After that the treasurers found the rest to be 421 kendinars and 3600 dahekans. Now, work out how much the whole treasure was worth.

    \noindent \textbf{Comments.} \textit{Marcian was Roman emperor of the East from 450 to 457. The kendinar and the dahekan were units of weight and currency. One kendinar is equal to 7200 dahekans.}

    \item The salaries of the clergy of St Sophia were distributed as follows. Deacons got the fifth share of it, priests got one tenth, bishops received two hundred litres and the rest of the monks got two thousand litres. Now, work out how many litres the total salary was. 

     \noindent \textbf{Comments.} \textit{St Sophia refers to Hagia Sophia of Constantinople. The litre was a unit of mass equivalent to one Roman pound.}

    \item The salaries of the cavalry are distributed as follows. One fifth was given to the honourable, one eighth to the seniors, and 150 kendinars were given to the rest of the cavalry. Now, work out how many kendinars the total salary was.

    \item I grew lettuce in my garden. A Roman entered for a walk and ate one fifth and one fifteenth of the lettuce. Because I knew of his gluttony, I drove him out and, counting the remaining lettuces, discovered that there were 110 left. Now, work out how many lettuces I had and how many the Roman ate.

    \item I was in Marmet, the princely estate of the Kamsarakans. I went to the bank of river Akhuryan  and saw a school of fish. With a cast net I caught a half, a quarter, and a seventh of the fish. The rest, 45 in number, were caught by a hand net. Now, work out how many fish there were in the school.

     \noindent \textbf{Comments.} \textit{The Akhuryan is a river in Armenia, a tributary of Aras (aka Yeraskh). Marmet was another name for the city of Yervandashat near the confluence of the Aras and Akhuryan rivers. See Figures \ref{fig:Sasanian Armenia} and \ref{fig:Ayrarat map}.
     }

    \item During the rebellion of the Armenians against the Persians, when Zorak Kamsarakan killed Suren, a member of the Armenian nobility sent a messenger to report this to the Persian King. The messenger walked 50 miles per day. Fifteen days later, when Zorak Kamsarakan learnt about that, he sent some people to pursue and apprehend him. They travelled 80 miles per day. Now, work out how many days later they would catch up with the messenger.

     \noindent \textbf{Comments.} \textit{See the comments on Problem \ref{prob:1}.}

    \item The Kamsarakans went on a hunting trip in Gen. They hunted many game animals and brought me a wild boar. It was huge, so I weighed it. It turned out that the entrails formed one quarter, the head was one tenth, the legs were one twentieth, and the tusks one ninetieth of the total weight. The body weighed 212 litres. Now, work out how much the whole wild boar weighed.

     \noindent \textbf{Comments.} \textit{Gen was in the canton of Chakatk in the Ayrarat province.}

    \item A sheatfish was caught in the Yeraskh river, near Marmet. I weighed it and found out that its head was one quarter of the total weight, and the tail was one sixth. The trunk weighed 140 litres. Now, work out how much the whole fish weighed.

    \item A merchant passed through three cities. In the first city a half and a third of his goods were levied as a duty. In the second city a half and a third of the remaining goods were levied, and similarly in the third city too. When he reached home, he had only 11 dahekans left. Now, work out how many dahekans he had at the beginning of his trip.

    \item I wanted to build a boat but I had nothing but three drams. I asked my relatives for help. They contributed respectively one third, one fourth, one sixth, one seventh, and one twenty eighth  of the total cost. Then I was able to construct the boat. Now, work out how much it cost.

    \item One of my students bought some excellent apples from the city of Khar and wanted to give them to me. He met three groups of jesters on his way. The first took a half and a quarter of all apples. The second took a half and a quarter of the remaining apples, and so did the third group too. There were five apples left after that which he gave to me. Now, work out how many apples my student had in total.

     \noindent \textbf{Comments.} \textit{Khar was an unidentified city in the Sasanian Empire.}

    \item There were a pot full of rose wine and three jars. I had the wine poured into the jars. The capacities of the jars were respectively equal to one third, one sixth, and one fourteenth of the capacity of the pot. At the end 54 jugs of wine were left in the pot. Now, work out how many jugs the total volume of the wine was.

    \item I had a thoroughbred horse. I sold it and used one quarter of the profit to buy cows, with one seventh I bought goats, then spent one tenth to buy oxen. The remaining 318 dahekans I used to buy sheep. Now, work out the price of the horse.

    \item I was building a church. I hired a bricklayer who could lay 140 bricks in one day. After 39 days I hired another bricklayer who could lay 218 bricks per day. When the number of bricks laid by the second bricklayer became equal to that of the first, the construction finished. Now, work out how many days it took for the number of bricks to become equal.

    \item A ship full of grain was sailing on the sea when a whale started chasing it. The sailors poured half of the grain into the sea as food for the whale out of fear. On the following three days they gave away respectively one fifth, one eighth, and one seventh of the rest of the grain. When they reached the port, they had only 7200 baskets of grain left. Now, work out how much grain there was in total.\footnote{Jonah would love this problem.}

    \item I had a metal pot which I broke into pieces and made other vessels. I used a third to make a cauldron, and a quarter to make another cauldron. I made two cups with a fifth of the pot, and two trays with a sixth. Finally, I used the remaining 210 drams to make a mug. Now, work out the weight of the metal pot.

    \item A man visited three churches. In the first church he told the God, ``Give me as much as I have now and I will give you 25 dahekans.'' He did the same in the second and third churches, after which he was left with nothing. Now, work out how many dahekans he had at the beginning.

    \item The hunting field of Nerseh Kamsarakan, the lord of Shirak and Arsharunik, was situated on the foot of mount Artin. One day it was invaded by herds of wild asses. The servants, having no hunting experience, ran to the village of Talin and reported this. Kamsarakan went there hastily with his brothers and some noblemen. They started hunting the asses. Half of the asses were trapped in snares, a quarter were shot with arrows, and the young, which formed one twelfth of the total, were captured alive. They also killed 360 asses with spears. Now, work out the total number of the asses. 

     \noindent \textbf{Comments.} \textit{Mount Artin and the Talin village were located in Shirak. Both are in the Aragatsotn region of present-day Armenia. See Figure \ref{fig:Ayrarat map}. Nerseh Kamsarakan lived in the 7th century and was a contemporary of Anania.}

    \item Arshavir's son Nerseh Kamsarakan, who was a namesake and ancestor of the above-mentioned Nerseh, waging war against the inhabitants of Bahl and defeating them, took many prisoners of war. When he reached the royal palace, he presented half of the captives to the Persian king, and one seventh of the remaining half to the king's son. Then he bid them farewell and returned to his country. On his way he visited the governor of the palace, who received him with great honours, not just like an aristocrat but like a king, and so he gave him one eighth of the remaining captives. Then he visited the chief military officer, who was called Khoravaran, and being honoured even more by him, he presented him one fourteenth of the remaining prisoners. Then he went on his way and eventually reached his country. There his younger brother, Hrahat, came to meet him. He gave Hrahat one thirteenth of the remaining prisoners. He gave one ninth of the rest of the captives to the aristocrats who came to greet him. Upon reaching Vagharshapat he gave one sixteenth of the remaining prisoners to the churches. Finally, he gave one twentieeth of the remaining captives to his elder brother, Sahak. As a result he kept only 570 prisoners for himself. Now, work out the total number of the prisoners.

     \noindent \textbf{Comments.} \textit{The sons of Arshavir Kamsarakan, Nerseh, Hrahat, and Sahak, took part in the rebellion of 481-484 against the Sasanian Empire led by Vahan Mamikonyan. Vagharshapat was (and still is) the spiritual capital of Armenia and the centre of the Armenian Apostolic Church.  See Figures \ref{fig:Sasanian Armenia} and \ref{fig:Ayrarat map}.}

    \item The pharaoh of Egypt was celebrating his birthday. He was in the habit of presenting 100 jugs of incense-flavoured wine to ten of his high ranking officials according to their merits and ranks. Now, divide the wine among the officials. \label{prob:22}
    
     \noindent \textbf{Comments.} \textit{It is presumed that the official ranked $k$-th among these ten gets $11-k$ times more than the official ranked 10th.}  

    \item There were 200 baskets of barley in a granary. Mice found their way into the granary and ate the whole barley. I caught and punished one of them. He confessed to having eaten only 80 grains. Now, work out how many grains of barley there were in the granary and how many mice ate it, \textbf{assuming that the capacity of one basket was 414720 grains}.

     \noindent \textbf{Comments.} \textit{I have added the last clause to make the problem complete and to make Anania's answer to this problem correct.}

    \item There was \textbf{a pool} in the city of Athens. Three pipes supplied water to the pool. The first could fill it up in 1 hour, the second in 2 hours, and the third in 3 hours. Now, work out in what fraction of an hour the three pipes would fill the pool together. \label{prob:24}

    \noindent \textbf{Comments.} \textit{In the original version three pools are mentioned at the beginning of the problem, but only one in the question. It is evident that there should be just one pool.}
\end{enumerate}

The following problems were discovered by Matevosyan in a parchment kept at Matenadaran \cite{matevosyan-norahayt-patarikner}. It contains eight complete and incomplete problems three of which are among Anania's problems (problems 17,18,24 above). Based on that and the similarity to Anania's problems the author concludes that those are also due to Anania. Below I present the three complete problems from \cite{matevosyan-norahayt-patarikner} which are not in the above list. To the best of my knowledge, these problems have not been translated into English.

\begin{enumerate}[resume]
    \item A sower from Arats sowed some seeds in his garden. A few days later three ants came and ate some of the seed. The first ate a third, the second a quarter and the third one twenty sixth of the seeds. After that the sower went to the garden and found 1239 seeds left. Now, work out the total number of the seeds.

    \item A man died and left his 500 dahekans to his three sons. The elder son got half of this, the middle and younger sons each got a quarter. Now, work out how much each of them received.

    \item A blind man passed by a clock and asked for the time. The answer was, ``Add to the current time half of one seventh of it and one quarter of it and that will complete the day.'' Now, work out what time it was when the blind man asked.
\end{enumerate}

\subsection{Solutions}

The following answers are in Anania's Book of Arithmetic. They are called ``solutions'' there and it is possible that originally complete solutions were provided of which only these answers are extant.

\begin{enumerate}
    \item Before the massacre, there were 1760 cavalry. 

    \item There were 144 pearls and their total value was 6720 drams.

    \item The treasure was 1686 kendinars.

    \item The total salary of the clergy was 3200 litres.

    \item The total salary of the cavalry was 240 kendinars. 
    
    \item There were 150 lettuces.

    \item There were 420 fish in total.

    \item They caught up in 25 days.

    \item The wild boar was 360 litres.

    \item The sheatfish was 240 litres.

    \item  The merchant had 2376 dahekans.

    \item The ship was 42 drams.

    \item There were 320 apples.

    \item There were 126 jugs of wine.

    \item The price of the horse was 616 dahekans.

    \item The bricklayers' bricks reached equality in 70 days.

    \item The bread in the ship was 24,000 baskets. 
    
    \item  The metal pot was 4200 dram.

    \item There were 21, 1/2, 1/4, 1/8 dahekans.

    \item There were 2160 asses in total.

    \item  There were 2240 prisoners in total.

    \item  The ten officials should respectively get the following shares of the total wine.
    
    1. 1, 1/2, 1/5, 1/10, 1/55

2. 3, 1/2, 1/10, 1/40, 1/88

3. 5, 1/3, 1/15, 1/44, 1/60, 1/66

4. 7, 1/5, 1/20, 1/44

5. 9, 1/11

6. 10, 1/2, 1/5, 1/10, 1/22, 1/30, 1/33

7. 12, 1/2, 1/10, 1/22, 1/30, 1/33, 1/55

8. 14, 1/3, 1/10, 1/15, 1/22

9. 16, 1/5, 1/10, 1/22, 1/55

10. 18, 1/12, 1/22, 1/33, 1/44

    \item  There were 8,294,4000 grains and 1,036,800 mice.\footnote{That was a really heavy infestation.}

    \item  The pipes would fill the pool in 1/4, 1/6, 1/12 and 1/22 of an hour.
\end{enumerate}

\subsection{Geography of Anania's problems}

\begin{figure}
    \centering
    \includegraphics[width=1\linewidth]{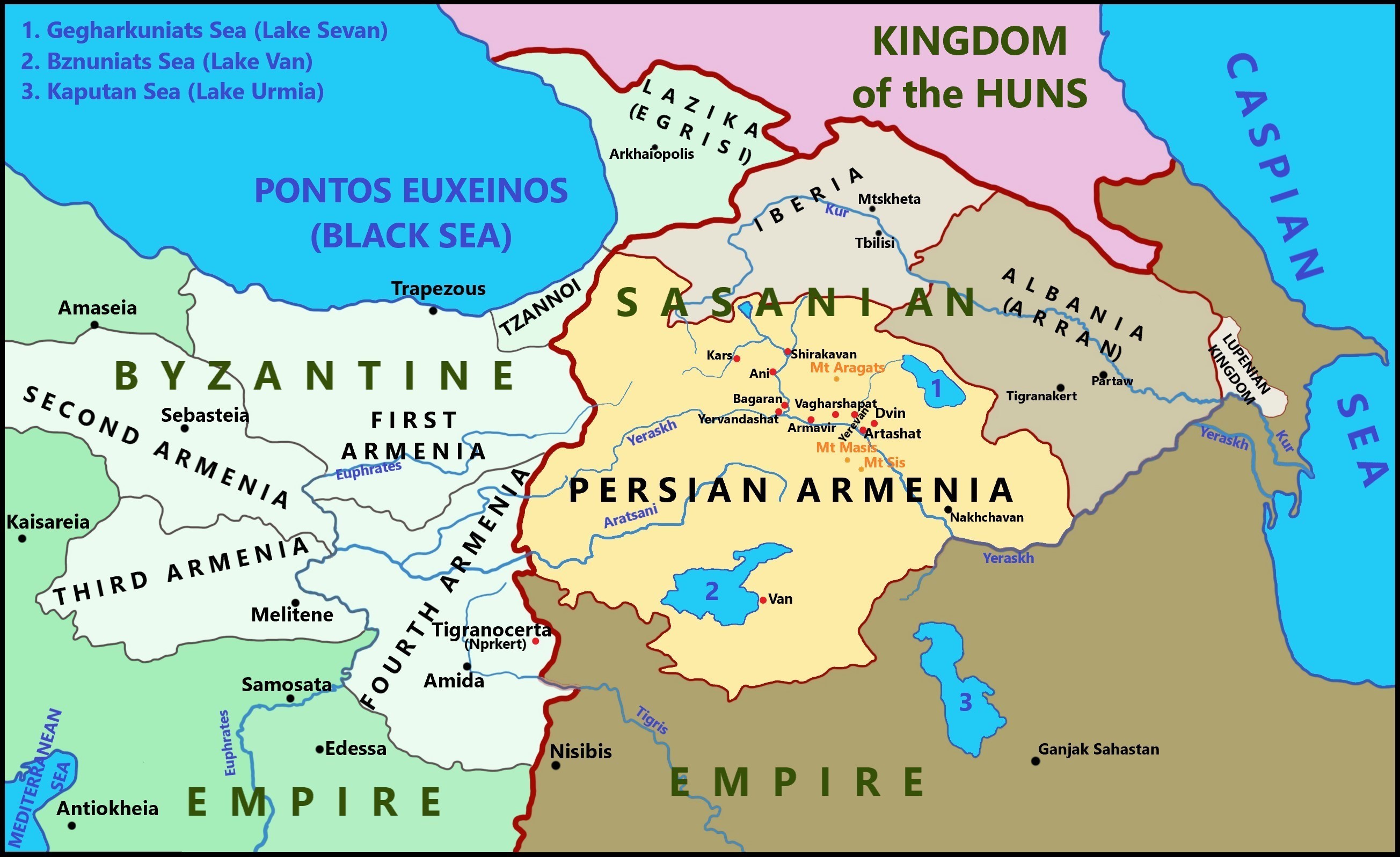}
    \caption{A map of Armenia and its surrounding areas between 387 and 591 AD. The red dots mark the 12 historical capitals of Armenia.}
    \label{fig:Sasanian Armenia}
\end{figure}

\begin{figure}
    \centering
    \includegraphics[width=1\linewidth]{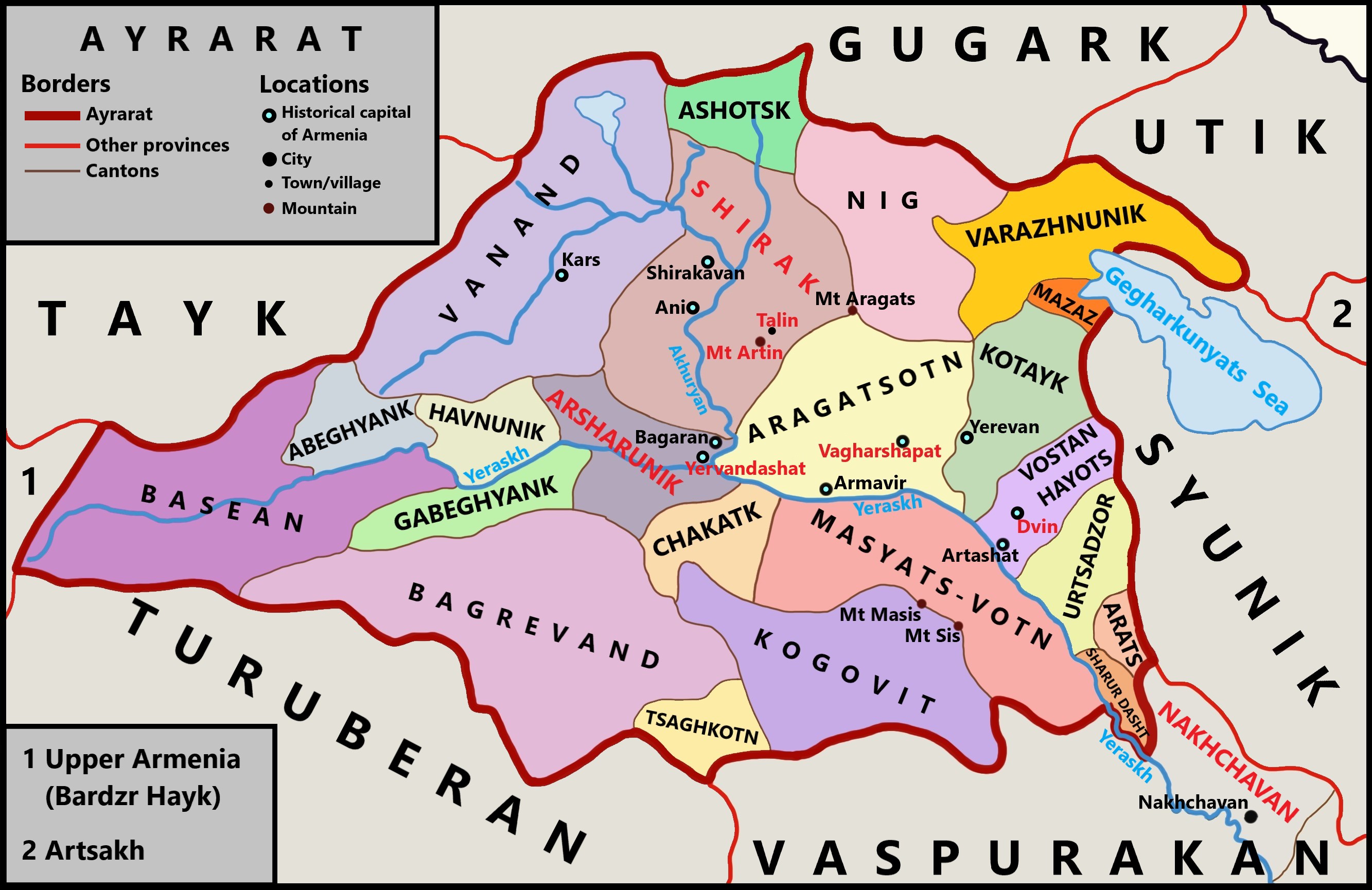}
    \caption{A map of Ayrarat and the surrounding provinces of Greater Armenia. Locations in red, and the Yeraskh and Akhuryan rivers, are mentioned in Anania's problems. }
    \label{fig:Ayrarat map}
\end{figure}

In this section I discuss the geography of Anania's problems. I made the maps in Figures \ref{fig:Sasanian Armenia} and \ref{fig:Ayrarat map} to highlight the locations mentioned in these problems. The first map shows Armenia after its partition between the Byzantine (Eastern Roman) and the Sasanian Empires in 387. The provinces Ayrarat, Turuberan, Syunik, Vaspurakan, Moxoene, Tayk of Greater Armenia formed Sasanian/Persian Armenia (Marzpanate of Armenia after 428) under Sasanian suzerainty. The provinces Sophene, Aghdznik, and Upper Armenia (Bardzr Hayk) became part of the Byzantine Empire. The other provinces -- Gugark, Utik, Artsakh, Paytakaran, Parskahayk, Corduene -- as well as Armenian Mesopotamia, were annexed by the Persians and included in the territories of the neighbouring political entities which, like Armenia, were under the suzerainty of the Sasanian Empire. Thus, Gugark was attached to Iberia, Artsakh and Utik were joined to Aghvank (aka Arran and Albania), while Paytakaran, Parskahayk, Corduene, and Armenian Mesopotamia were under the direct control of the Sasanian Empire.

Byzantine Armenia comprised the three above-mentioned provinces (Sophene, Aghdznik, Upper Armenia) along with Lesser Armenia. It was divided into four administrative regions -- First Armenia, Second Armenia, Third Armenia, and Fourth Armenia. Anania's teacher Tychicus was based in Trebizond (Trapezous) in First Armenia, where Anania spent 8 years mastering the basics of arithmetic.

This map significantly changed in 591 after the second partition of Armenia when most of its territory was annexed by the Byzantine Empire. Anania lived after this, in the 7th century, however the historical events mentioned in his problems took place in Sasanian Armenia and the Sasanian Empire between 387 and 591 AD.

The second map shows the Ayrarat province of Greater Armenia, along with its neighbouring provinces (partially). Anania was born in the Shirak canton of Ayrarat (hence he was called \textit{Shirakatsi}, i.e. \textit{of Shirak}), and almost all of the locations mentioned in the problems are in Ayrarat and are highlighted in this map. Nakhchavan, mentioned in Problem 1, was a canton of the Vaspurakan province and bordered the Arats and Sharur Dasht cantons of Ayrarat. Note that the city of Yervandashat in Arsharunik was also known as Marmet, which is the name used by Anania.

The maps also show a few important places that do not appear in the problems, namely Mounts Masis, Sis, and Aragats, and the twelve historical capitals of Armenia (ten of which were located in Ayrarat).

I used three sources when creating these maps. The first is based on the map ``The first Byzantine expansion into Armenia, 387-591'' from \textit{Armenia: A Historical Atlas} by Robert Hewsen \cite[map 65]{Hewsen-historical-atlas}. The second is based on the map of Ayrarat in the \textit{Armenian Soviet Encyclopedia} \cite[volume 1, page 353]{ASE} and the map of Greater Armenia in Suren Yeremyan's \textit{Armenia according to ``Ashkharhatsuyts''} \cite[pp. 42--43]{yeremyan-armenia-ashkhahatsuyts}. All of these maps are in turn based on \textit{Ashkharhatsuyts}, Anania Shirakatsi's \textit{Geography} \cite{Hewsen-ashkharhatsuyts}. It is one of Anania's most important works and contains detailed information about Greater Armenia, Persia, and the Caucasus (but no actual maps). 

\subsection{Analysis of Anania's problems and solutions} 

One important questions is whether Anania composed these problems himself or translated them from another source. There are several reasons to believe that the problems are mostly due to Anania. Firstly, Matthews \cite{matthews-anania-of-shirak} writes that the problems and solutions and the entertaining puzzles (see the next section) ``are possibly the oldest extant texts of their kind in the world''. If no earlier text with these kind of problems is known, Anania's problems are unlikely to be translated or adapted from a written source.

Secondly, most of the problems are about Armenia, Armenians and their everyday life, some are based on historical events, and others on units of measurement used in Armenia in the early Middle Ages. It seems that the numbers and amounts appearing in the problems are chosen to make the narrative as plausible as possible in the context of early medieval Armenia. It is of course quite possible that he studied similar problems with his teacher Tychicus and adapted or even copied some of his problems (e.g. \ref{prob:22} and \ref{prob:24}) from these, but even in that case it seems more likely that he would have to change the numbers appearing in these problems to match their context (e.g. using Armenian units would make this necessary if one wants to get whole number solutions, see below). Therefore, most of the problems were most likely composed by Anania himself and the solutions were also due to him.

Thirdly, in the text containing the arithmetic tables, Anania gives a short introduction before the tables of each arithmetic operation, and at the very beginning he writes that he is presenting the arithmetic developed by his ``ancestors'', probably referring to Tychicus and other mathematicians whose work he was familiar with. Thus, he clearly acknowledges that he is presenting something developed by others. There seems to be no such statement in the manuscript containing the problems and solutions. More precisely, Anania's work has survived in various copies of his original texts, and often different copies of the same text disagree at certain places. In particular, Greenwood \cite{greenwood-reassessment} writes that one of the manuscripts (Matenadaran 3078) does give an introduction before the problems and, in particular, states, ``Through making a demonstration, I shall try [to supply this] by means of a few [examples] in the manner of an introduction through a treatise, \textbf{passing over errors} for the sake of coinciding properly with the most profound.'' I do not have access to the original manuscript but I have not seen such an introduction in any other paper or book on Anania's work, which are presumably based on other manuscripts. Actually, given Greenwood's translation, it is quite suspicious that Anania, writing a textbook for his students, would decide to keep the errors of the alleged source(s). Therefore, I find it more plausible that this introduction was written by the author of that particular manuscript who copied from Anania's original text or possibly from another copy. Then ``the most profound'' would refer to Anania, and it would make sense for a copier to choose to keep Anania's errors for historical reasons.

Now let us discuss the mathematical contents of Anania's problems. Both the problems and solutions indicate that Anania worked solely with Egyptian fractions and not with other forms of rational numbers. An Egyptian fraction is a finite sum of reciprocals of distinct positive integers. For instance, instead of writing three quarters Anania writes a half and a quarter. It is indeed true that every positive rational number can be written as an Egyptian fraction. Historically, different methods have been used by different authors to achieve this  \cite{knorr-egypt-greece}, but it is not clear how Anania represented a positive rational number as an Egyptian fraction.

All twenty four problems of Anania can be solved by a linear equation in a single variable. So it is clear that he possessed the knowledge for solving such equations over the positive rational numbers. What is more interesting though is that this also implies that Anania was able to manipulate Egyptian fractions and, in particular, write the sum, product, difference, and ratio of two Egyptian fractions as Egyptian fractions. This is quite more advanced than working with other representations of rational numbers and requires some non-trivial arithmetic skills. Thus, we conclude that Anania possessed these skills, even though he did not write about them (or his writings did not survive). Moreover, the answers to most questions are whole numbers, and in many cases they must be, as they describe the number of people or animals or other things. Therefore, the numbers used in the problems should have been chosen in a way that would result in integer solutions. This also testifies to Anania's advanced arithmetic skills.

Thus, we conclude that while Anania's Book of Arithmetic is much less elaborate than the works of ancient mathematicians (e.g. Euclid) and his contemporaries (e.g. Brahmagupta), his proficiency in arithmetic was still quite sophisticated.

\section{Entertaining puzzles}\label{sec: entertaining puzzles}

The third section of Anania's Book of Arithmetic contains several entertaining puzzles which people could use during feasts to challenge each other's mathematical skills. Their purpose was to raise public interest in arithmetic. The Armenian word that Anania used to describe these puzzles is ``\textarmenian{Խրախճանականք}'' (\textit{Khrakhchanakank}), which can be translated as ``feast riddles'' or ``feast puzzles''. Petrosyan studied the cultural aspects of these riddles in \cite{petrosyan-xraxhcanakanq}, where he also presented English translations of two of them (2 and 5 in the below list). To the best of my knowledge, the other puzzles have never appeared in English before.

Here I present the full list of Anania's puzzles translated from \cite{abrahamyan-petrosyan-shirakaci-matenagrutyun}. In fact, there are three more puzzles, but these are either incomplete or have a somewhat different nature, and so are omitted from the below list.

\begin{enumerate}[itemsep=1ex]
    \item Tell your friend, ``I can work out when you want to have dinner and how much wine you want to drink.'' If he agrees to play, tell him, ``Keep the number of hours when you want to have dinner in your mind. Double it, then add 5. Then multiply by 5, add 10, and multiply the result by 10. Then add the number of glasses of wine you want to drink.'' Then ask him for the resulting number. Whatever number he says, subtract 350. The number of hundreds in this new number is the hour when he wants to dine. The remainder (upon division by 100) is the number of glasses he wants to drink. If your friend is unskilled and the number of glasses of wine turns out to be 100 or more, tell him that it is not possible to drink 100 glasses of wine in one hour. 

    \item Tell your friend that once upon a time during a feast a Persian tourist noticed a group of Greek tourists, called them and said, ``If you join me, then another group equal to yours in number joins us, then one more group half your size, then one more group one quarter of your number, and counting me too, we will be 100 people.'' Then have your friend work out the number of Greek tourists. If you friend is smart, he will quickly find the number of the Greeks to be 36. Otherwise you can make fun of his struggles trying to work out the number.

    \item Tell your friend that you can find out how much money he has in his pouch. When he says, ``Go ahead,'' tell him, ``Pick the quantity of your dram, add an equal amount to it, then double the sum, then add the original amount to it, and double it again.'' When he tells you the number he obtained, regardless of whether it is odd or even, divide it by ten and you will get the amount of money he had in his pouch.

    \item Tell your friend that a Hun looked after your chickens for 100 years and ate 100 eggs each day.\footnote{A hundred eggs a day! It seems Anania didn't like Huns.} Then ask him to work out the total number of eggs. If he is smart, he will tell you that the total number is 365 myriads\footnote{\textit{Myriad} (\textit{\textarmenian{բյուր}} (\textit{byur}) in Armenian) means 10,000.}. If he is foolish, his suffering will cause you great joy.

    \item Tell your friend, ``If you sell 60 skins of wine at 2 drams for every 5 skins, how much money will you get?'' He will answer, ``24 drams.'' Then say to him that you can sell the wine at the same price and get one more dram. This is how; you will divide 60 by two, then you will divide the first 30 skins into groups of three and sell each group for 1 dram. Then you will divide the remaining 30 into groups of two and sell each for 1 dram. Together, the price is still 1 dram for 5 skins, but you will get 25 drams, thus shaming him and rejoicing at his bafflement.
\end{enumerate}

These puzzles support the point from the previous section about Anania's sophisticated arithmetic skills. The first puzzle demonstrates that apart from masterfully manipulating natural and rational numbers, Anania was also familiar with identities (i.e. he could show that an equality holds for all possible values of the variables involved), as well as division with remainder. The second one is like most of his problems, with Egyptian fractions, but presumably it was considered easy enough by Anania that it could be solved impromptu during a feast. The third puzzle is about a linear identity, while the fourth presents a simple multiplication exercise. The fifth puzzle is a trick question and is arguably the most interesting and entertaining one. I offered it to some students, and even the most capable ones among them had to think for a few minutes to understand where the trick was. Some of the others are probably still baffled by it...

\bibliographystyle{alpha}
\bibliography{refhist}

\end{document}